\documentclass{article}
\usepackage{amsmath}
\usepackage{amssymb}
\usepackage{amsthm}
\usepackage{comment}

\usepackage{xcolor}

\newcommand{\F}{\mathbb{F}}
\newcommand{\Z}{\mathbb{Z}}

\newtheorem{theorem}{Theorem}
\newtheorem{lemma}{Lemma}

\newtheorem{corollary}{Corollary}

\title{Normality of the Thue-Morse function for finite fields along polynomial values}
\author{Mehdi Makhul and Arne Winterhof}
\date{}

\begin{document}
	
	\maketitle
	
	\begin{center}
		Johann Radon Institute for Computational and Applied Mathematics,
		Austrian Academy of Sciences, Altenberger Str. 69, 4040 Linz, Austria\\
		E-mail: \{mehdi.makhul,arne.winterhof\}@oeaw.ac.at
	\end{center}
	
	\abstract{Let $\F_q$ be the finite field of $q$ elements, where $q=p^r$ is a power of the prime $p$, and $\left(\beta_1, \beta_2, \dots, \beta_r \right)$ be an ordered basis of $\F_q$ over $\F_p$. For $$\xi=\sum_{i=1}^rx_i\beta_i, \quad \quad x_i\in\F_p,$$ we define the Thue-Morse or sum-of-digits function $T(\xi)$ on $\F_q$ by
		\[
		T(\xi)=\sum_{i=1}^{r}x_i.
		\]
		For a given pattern length $s$ with $1\le s\le  q$, a subset ${\cal A}=\{\alpha_1,\ldots,\alpha_s\}\subset \F_q$, a polynomial $f(X)\in\F_q[X]$ of degree $d$ and a vector $\underline{c}=(c_1,\ldots,c_s)\in\F_p^s$
		we put 
		\[
		{\cal T}(\underline{c},{\cal A},f)=\{\xi\in\F_q : T(f(\xi+\alpha_i))=c_i,~i=1,\ldots,s\}.
		\]
		In this paper we will see that  under some natural conditions, the size of~${\cal T}(\underline{c},{\cal A},f)$ is asymptotically the same for all~$\underline{c}$ and ${\cal A}$ in both cases, $p\rightarrow \infty$ and $r\rightarrow \infty$, respectively. More precisely, we have 
		\[
		\left||{\cal T}(\underline{c},{\cal A},f)|-p^{r-s}\right|\le (d-1)q^{1/2}\]
		under certain conditions on $d,q$ and $s$.
		For monomials of large degree we improve this bound as well as we find conditions on $d,q$ and $s$ for which this bound is not true. In particular, if $1\le d<p$ we have the dichotomy that the bound is valid if $s\le d$ and fails for some $\underline{c}$ and ${\cal A}$ if $s\ge d+1$. The case $s=1$ was studied before by Dartyge and S\'ark\"ozy.
	}\\
	
	
	
	\section{Introduction}
	
	\subsection{The problem for binary sequences}
	
	For positive integers $M$ and $s$, a binary sequence $(a_n)$
	and a binary pattern 
	$${\cal E}_s=(\varepsilon_0,\ldots,\varepsilon_{s-1})\in \{0,1\}^s$$
	of length $s$ we denote by $N(a_n,M,{\cal E}_s)$ the number of $n$ with $0\le n<M$
	and $(a_n,a_{n+1},\ldots,a_{n+s-1})={\cal E}_s.$
	The sequence $(a_n)$ is {\em normal} if for any fixed $s$ and any pattern ${\cal E}_s$
	of length $s$,
	$$\lim_{M\rightarrow \infty}\frac{N(a_n,M,{\cal E}_s)}{M}=\frac{1}{2^s}.$$
	
	The {\em Thue-Morse} or {\em sum-of-digits sequence} $(t_n)$ is defined by
	$$t_n=\sum_{i=0}^\infty n_i \bmod 2,\quad n=0,1,\ldots$$
	if
	$$n=\sum_{i=0}^\infty n_i2^i,\quad n_0,n_1,\ldots \in \{0,1\},$$
	is the binary expansion of $n$.
	Recently, Drmota, Mauduit and Rivat \cite{DMR19} showed that the Thue-Morse sequence along squares, that is, $(t_{n^2})$ is normal. 
	It is conjectured but not proved yet that the subsequence of the Thue-Morse sequence along any polynomial of degree $d\ge 3$ is normal as well, see \cite[Conjecture 1]{DMR19}.
	Even the weaker problem of determining the frequency of $0$ and $1$ in the subsequence of the Thue-Morse sequence along any polynomial of degree $d\ge 3$ seems to be out of reach, see 
	\cite[above Conjecture 1]{DMR19}.
	
	However, the analog of the latter weaker problem for the Thue-Morse sequence in the finite field setting was settled by Dartyge and S\'ark\"ozy \cite{DS13}.
	
	\subsection{The analog for finite fields}
	
	This paper deals with the following analog of the normality problem. Let~$q=p^r$ be the power of a prime $p$ 
	and 
	$${\cal B}=( \beta_1,\ldots ,\beta_r)$$
	be an ordered basis of the finite field $\F_q$ over $\F_p$.
	Then any $\xi\in\F_q$ has a unique representation  
	$$
	\xi=\sum_{j=1}^r x_j\beta_j \quad\mbox{with } x_j\in \F_p,\quad j=1,\ldots,r.
	$$
	The coefficients $x_1,\ldots , x_r$ are called the {\em digits} with respect to the basis~${\cal B}$.
	
	Dartyge and S\'ark\"ozy \cite{DS13} introduced the {\em Thue-Morse} or {\em sum-of-digits function} $T(\xi)$ for the finite field $\F_q$ with respect to the basis ${\cal B}$:
	$$
	T(\xi)=\sum_{i=1}^{r}x_i,\quad \xi=x_1\beta_1+\cdots +x_r\beta_r\in \F_q.
	$$
	For a given pattern length $s$ with $1\le s\le  q$, a subset 
	$${\cal A}=\{\alpha_1,\ldots,\alpha_s\}\subset \F_q,$$ 
	a polynomial $f(X)\in\F_q[X]$ 
	and a vector $\underline{c}=(c_1,\ldots,c_s)\in\F_p^s$
	we put 
	$$
	{\cal T}(\underline{c},{\cal A},f)=\{\xi\in\F_q : T(f(\xi+\alpha_i))=c_i,~i=1,\ldots,s\}.
	$$
	In~\cite{DS13} the Weil bound, see Lemma~\ref{weil} below, was used to bound the cardinality of ${\cal T}(\underline{c},{\cal A},f)$ for~$s=1$:
	
	Let $f(X)\in \F_q[X]$ be a polynomial of degree $d$. Then for all $c \in \F_p$
	\begin{equation}\label{DarSar}
	\left||{\cal T}(c,f)|-p^{r-1}\right|\le (d-1)q^{1/2},\quad  \gcd(d,p)=1,
	\end{equation}
	where
	\[
	{\cal T}(c,f)=\{\xi\in\F_q : T(f(\xi))=c\}.
	\]
	Note that the condition $\gcd(d,p)=1$ can be relaxed to the condition that~$f(X)$ is not of the form 
	$g(X)^p-g(X)+c$ for some $g(X)\in \F_q[X]$ and $c\in \F_q$. For example, $f(X)=X^p$ is not of the form $g(X)^p-g(X)+c$ but does not satisfy $\gcd(d,p)=1$.
	
	Our goal is to prove that, under some natural conditions, the size of~${\cal T}(\underline{c},{\cal A},f)$ is asymptotically the same for all~$\underline{c}$ and ${\cal A}$.
	
	\subsection{Results of this paper}
	
	First we study monomials and prove the following result in Section~\ref{monomials}.
	
	\begin{theorem}\label{Xd}
		Let $d$ 
		be any integer with $1\le d<q$ with unique representation
		$$d=(d_0+d_1p+\ldots+d_{n-1}p^{n-1})\gcd(d,q)$$ 
		where
		$$1  \le n\le r-\frac{\log(\gcd(d,q))}{\log p}, \quad 0\le d_i<p,~i=0,\ldots,n-1,~d_0d_{n-1}\ne 0.$$  
		Let denote by
		$$f_d(X)=X^d\in \F_q[X]$$ 
		the monomial of degree $d$ and if $n\ge 2$, assume 
		$$d_m=d_{m+1}=\ldots=d_{m+k-1}=p-1$$
		for some $m$ and $k$ with 
		$$1\le m\le m+k\le n-1.$$
		
		1. For any positive integer
		\begin{equation}\label{scond} s\le\left\{\begin{array}{cc}\max\left\{d_0,(d_{m+k}+1)p^k\right\}, & n\ge 2,\\
		d_0, & n=1,\end{array}\right.
		\end{equation}
		any ${\cal A}\subset \F_q$ of size $s$ and any $\underline{c}\in \F_p^s$ we have
		$$\left||{\cal T}(\underline{c},{\cal A},f_d)|-p^{r-s}\right|\le \left(\frac{d}{\gcd(d,q)}-1\right)q^{1/2}.$$
		
		2. Conversely, if
		\begin{equation}\label{Dcond}(d_0+1)(d_1+1)\cdots (d_{n-1}+1)\le p,
		\end{equation}
		for any integer $s$ with
		\begin{equation}\label{scond2} q\ge s\ge (d_0+1)(d_1+1)\cdots (d_{n-1}+1),
		\end{equation}
		there is  ${\cal A}\subset \F_q$ of size $s$ and a vector $\underline{c}\in \F_p^s$ for which ${\cal T}(\underline{c},{\cal A},f_d)$
		is empty. 
		
		3. For any $s$ with
		\begin{equation}\label{Dcond2} q\ge s>((d_0+1)(d_1+1)\cdots (d_{n-1}+1)-1)r
		\end{equation}
		and any ${\cal A}\subset \F_q$ of size $s$
		there is a vector $\underline{c}\in \F_p^s$ for which ${\cal T}(\underline{c},{\cal A},f_d)$
		is empty.
	\end{theorem}
	
	For $d<p$ we have the following dichotomy:
	\begin{corollary} Assume $1\le d<p$.
		For $s\le d$ we have 
		$$\left||{\cal T}(\underline{c},{\cal A},f_d)|-p^{r-s}\right|\le (d-1)q^{1/2}$$
		and for any $s$ with $q\ge s>d$   
		there is ${\cal A}\subset \F_q$ of size $s$ and a vector $\underline{c}\in \F_p^s$ for which ${\cal T}(\underline{c},{\cal A},f_d)$
		is empty. 
	\end{corollary}

	Theorem~\ref{Xd} provides two asymptotic formulas for $|{\cal T}(\underline{c},{\cal A},X^d)|$
	for $r\rightarrow \infty$ and $p\rightarrow \infty$, respectively. 
	
	Assume that $p$, $j$, $n$, $d=(d_0+d_1p+\ldots+d_{n-1} p^{n-1})p^j$ and $s$
	satisfying~$(\ref{scond})$ are fixed. Then we have
	$$\lim_{r\rightarrow \infty}\frac{|{\cal T}(\underline{c},{\cal A},f_d)|}{p^{r-s}}=1$$
	for any $\underline{c}\in \F_p^s$ and ${\cal A}\subset \F_q$ of size $|{\cal A}|=s$.
	We may say that $T(f_d)$ is~{\em $r$-normal} if $(\ref{scond})$ is satisfied.
	
	Assume that $j=0$ and $d$, $r$ and $s$ are fixed with $1\le s\le \min\{d,\lfloor (r-1)/2\rfloor\}$.
	Then we have
	$$\lim_{p\rightarrow \infty}\frac{|{\cal T}(\underline{c},{\cal A},f_d)|}{p^{r-s}}=1$$
	for any $\underline{c}\in \F_p^s$ and ${\cal A}\subset \F_q$ of size $|{\cal A}|=s$.
	We may say that $T(f_d)$ is~{\em $p$-normal} for $1\le s\le \min\{d,\lfloor (r-1)/2\rfloor\}$.

	Theorem~\ref{Xd} is only non-trivial for small degrees. 
	However, for very large degrees we prove the following non-trivial result in Section~\ref{sec:rational}. 
	
	\begin{theorem}\label{inv}
		Let $f_{q-1-d}(X)=X^{q-1-d}$ be a monomial of degree $q-1-d$ with $1\le d< q$. Then for any ${\cal A}\subset \F_q$ of size $s$ and any $\underline{c}\in \F_p^s$, we have
		\[
		\left||{\cal T}(\underline{c},{\cal A},f_{q-1-d})|-p^{r-s}\right|\le \left(\left(\frac{d}{\gcd(d,q)}+1\right)s-2\right)q^{1/2}+s+1. 
		\]
	\end{theorem}
	Note that with the convention $0^{-1}=0$ we have 
	$$\xi^{q-1-d}=\xi^{-d}\quad \mbox{for any }\xi\in \F_q$$ 
	and can identify the monomial $f_{q-1-d}(X)=X^{q-1-d}$ with the rational function $f_{-d}(X)=X^{-d}$. However, the latter representation is independent of $q$ and we can  state two asymptotic formulas for $|{\cal T}(\underline{c},{\cal A},f_{-d})|$ as well.
	
	For any fixed $d$, $p$ and $s$ we have
	$$\lim_{r\rightarrow \infty}\frac{|{\cal T}(\underline{c},{\cal A},f_{-d})|}{p^{r-s}}=1,$$
	that is, $T(f_{-d})$ is $r$-normal.
	
	For any fixed $d$, $s$ and $r$  with $1\le s\le \lfloor (r-1)/2\rfloor$ we have
	$$\lim_{p\rightarrow \infty}\frac{|{\cal T}(\underline{c},{\cal A},f_{-d})|}{p^{r-s}}=1,$$
	that is, $T(f_{-d})$ is $p$-normal for $1\le s\le \lfloor (r-1)/2\rfloor$.
	
	Finally, we extend our results to arbitrary polynomials in Section~\ref{secarb}. 
	\begin{theorem}\label{arb}
		Let $d$ 
		be any integer with $1\le d<q$ and $\gcd(d,q)=1$. Denote 
		$$d_0\equiv d\bmod p,\quad 1\le d_0<p.$$ 
		Let $f(X)\in \F_q[X]$ be any polynomial of degree $d$.
		
		1. For any integer $s$ with
		$$1\le s\le d_0,$$
		any ${\cal A}\subset \F_q$ of size $s$ and any $\underline{c}\in \F_p^s$ we have
		$$\left||{\cal T}(\underline{c},{\cal A},f)|-p^{r-s}\right|\le (d-1)q^{1/2}.$$
		
		2. Conversely, if $f(X)\in \F_p[X]$,
		then
		for any integer $s$ with
		$$q\ge s\ge d+1,$$
		there is  ${\cal A}\subset \F_q$ of size $s$ and a vector $\underline{c}\in \F_p^s$ for which ${\cal T}(\underline{c},{\cal A},f)$
		is empty. 
		
		3. For any $f(X)\in \F_q[X]$, any $s$ with
		$$q\ge s\ge dr$$
		and any ${\cal A}\subset \F_q$ of size $s$
		there is a vector $\underline{c}\in \F_p^s$ for which ${\cal T}(\underline{c},{\cal A},f)$
		is empty.
	\end{theorem}
	The restriction $\gcd(d,p)=1$ is needed. We give counter-examples of degree $d$ with $\gcd(d,p)>1$ in Section~\ref{counterp}.
	
	Again, for $f(X)\in \F_p[X]$ and $1\le d<p$ we have a dichotomy.

	Moreover, for any fixed $d$, $p$ and $s$ with $\gcd(d,q)=1$ and $1\le s\le d_0$ and any $f(X)\in \F_p[X]$ of degree $d$, $T(f)$ is $r$-normal. Note that any~$f(X)\in \F_p[X]$ is an element of $\F_{p^{r}}[X]$ for $r=1,2,\ldots$
	
	For fixed $d$, $r$ and $s$ with $1\le s\le \min\{d,\lfloor (r-1)/2\rfloor\}$ and any $f(X)\in \Z[X]$ of degree $d$, $T(f)$ is $p$-normal. Here $f(X)\in \Z[X]$ can be identified with an element of $\F_p[X]$ for all primes $p$.

	We start with a section on preliminary results used in the proofs. Then we show that
	\begin{equation}\label{main} 
	\left||{\cal T}(\underline{c},{\cal A},f)|-p^{r-s}\right|\le (\deg(f)-1)q^{1/2}
	\end{equation}
	under certain conditions in Section~\ref{tretal}. 
	In Sections~\ref{monomials} to \ref{secarb} we show that these conditions are fulfilled under the assumptions of our theorems. We finish the paper with some remarks on related work in Section~\ref{final}.

	\section{Preliminary results}
	
	We start with the Weil bound, see \cite[Theorem 5.38 and comments below]{LN97}, \cite[Theorem 2E]{sc} or~\cite{we48}.
	\begin{lemma}\label{weil}
		Let $\psi$ be the additive canonical character of the finite field~$\F_q$, and $f(X)$ be a polynomial of degree $d\ge 1$ over $\F_q$, which is not of the form $g(X)^p-g(X)+c$ for some polynomial $g(X)\in \F_q[X]$ and $c\in \F_q$. Then we have
		\[
		\left|\sum_{\xi \in \F_q}\psi\left(f(\xi)\right)\right| \le (d-1)q^{1/2}.
		\]
	\end{lemma}
	
	We also use the analog of the Weil bound for rational functions 
	$$\frac{f(X)}{g(X)}\in \F_q(X)$$ 
	of C.\ Moreno and O.\ Moreno~\cite[Theorem $2$]{momo}.
	We only need the special case that $\deg(f)\le \deg(g)$. 
	
	
	\begin{lemma}\label{th:momo} 
		Let $\psi$ be a nontrivial additive character of $\F_q$ and let $\frac{f(X)}{g(X)}\in \F_q(X)$ be a rational function over $\F_q$. Let $s$ be the number of distinct roots of the polynomial $g(X)$ in the algebraic closure $\overline{\F_q}$ of $\F_q$. Suppose that~$\frac{f(X)}{g(X)}$ is not of the form $H(X)^p-H(X)$, where $H(X)$ is a rational function over~$\overline{\F_q}$. If $\deg(f)\le \deg(g)$, then we have
		\[
		\left|\sum_{\xi \in \F_q, g(\xi)\not=0}\psi\left(\frac{f(\xi)}{g(\xi)}\right)\right| \le  (\deg(g)+s-2)\sqrt{q}+1.
		\]
	\end{lemma}
	Note that $g(X)^p-g(X)+c$ with $g(X)\in \F_q(X)$ and $c\in \F_q$ can be written as 
	$h(X)^p-h(X)$ for $h(X)=g(X)+\gamma\in \overline{\F_q}(X)$, where~$\gamma\in \overline{\F_q}$ is a zero of the polynomial $X^p-X-c$.\\
	
	Next we state Lucas' congruence, see \cite{EL} or \cite[Lemma~6.3.10]{niwi15}.
	\begin{lemma}\label{lm:lucas}
		Let $p$ be a prime.
		If $m$ and $n$ are two natural numbers with $p$-adic expansions
		\[
		m=m_{r-1}p^{r-1}+m_{r-2}p^{r-2}+\ldots+m_1p+m_0,\quad 0 \le m_0,\ldots,m_{r-1}< p,\]
		and
		\[n=n_{r-1}p^{r-1}+n_{r-2}p^{r-2}+\ldots+ n_1p+n_0,\quad 0 \le n_0,\ldots,n_{r-1}< p,
		\]
		then we have
		\[
		{m \choose n} \equiv \prod_{j=0}^{r-1} {m_j \choose n_j} \mod p.
		\]
	\end{lemma}
	
	As a consequence of Lucas' congruence we can count the number of nonzero binomials coefficients ${m\choose n }\bmod p$
	for fixed $m$. Indeed, by Lucas' congruence 
	$${m \choose n} \not\equiv 0 \mod p \mbox{ if and only if } {m_j \choose n_j} \not\equiv 0 \mod p \mbox{ for $j=0,\ldots,r-1$},$$ or equivalently, 
	$$0\le n_j\le m_j\quad \mbox{for }j=0,\ldots,r-1.$$
	Therefore, we have the following result of Fine \cite[Theorem $2$]{NF}:
	
	\begin{lemma}\label{th:fine} Let $p$ be a prime and $m$ an integer with $p$-adic expansion $$m=m_{r-1}p^{r-1}+m_{r-2}p^{r-2}+\ldots+ m_1p+m_0, \quad 0\le m_0,\ldots,m_{r-1}<p.$$ Then the number of nonzero binomial coefficients ${m \choose n}\bmod p$ with $0\le n\le m$ is
		\[
		\prod_{j=0}^{r-1}(m_j+1).
		\]
	\end{lemma}

	\section{Trace, dual basis and exponential sums}\label{tretal}
	
	Let
	$$
	{\rm Tr}( \xi)=\sum_{i=0}^{r-1}\xi^{p^i}\in\F_p
	$$ 
	denote the (absolute) \emph{trace} of $\xi\in \F_q$. 
	Let $(\delta_1,\ldots,\delta_r)$ denote the (existent and unique) \emph{dual basis} of the basis~${\cal B} =( \beta_1,\ldots ,\beta_r)$ of $\F_q$, see for example~\cite{LN97}, 
	that is,  
	\begin{equation*}
	{\rm Tr}(\delta_i\beta_j)=
	\begin{cases}
	1 & {\rm if}\ i=j,\\
	0 & {\rm if}\ i\not =j, 
	\end{cases}
	\quad 1\le i,j\le r.
	\end{equation*}
	Then we have 
	$${\rm Tr}(\delta_i\xi)=x_i\quad\mbox{for any}\quad \xi=\sum_{j=1}^r x_j \beta_j\in\F_q\quad \mbox{with }x_j\in \F_p,
	$$
	and
	$$T(\xi)={\rm Tr}(\delta \xi),\quad \mbox{where }\delta=\sum_{i=1}^r\delta_i.$$
	Note that 
	$$\delta\ne 0$$ 
	since $\delta_1,\ldots,\delta_r$ are linearly independent.
	
	Put 
	$$e_p(x)=\exp\left(\frac{2\pi i x}{p}\right)\quad \mbox{for }x\in \F_p.$$
	Since
	$$\sum_{a\in \F_p}e_p(ax)=\left\{\begin{array}{cc} 0,& x\ne 0,\\
	p, & x=0,\end{array}\right. \quad x\in \F_p,$$
	we get
	\begin{eqnarray*} |{\cal T}(\underline{c},{\cal A},f)|&=&
		\frac{1}{p^s} 
		\sum_{\xi\in \F_q}\prod_{i=1}^s\sum_{a\in \F_p}e_p\left(a (T(f(\xi+\alpha_i))-c_i)\right)\\
		&=&\frac{1}{p^s} \sum_{a_1,\ldots,a_s\in \F_p}
		\sum_{\xi\in \F_q}e_p\left(\sum_{i=1}^sa_i (T(f(\xi+\alpha_i))-c_i)\right).
	\end{eqnarray*}
	Separating the term for $a_1=\ldots=a_s=0$ we get 
	\begin{equation}\label{redchar}\left||{\cal T}(\underline{c},{\cal A},f)|-p^{r-s}\right|\le \max_{(a_1,\ldots,a_s)\ne (0,\ldots,0)} \left|\sum_{\xi\in \F_q} \psi(F_{a_1,\ldots,a_s}(\xi))\right|,
	\end{equation}
	where
	$$\psi(\xi)=e_p({\rm Tr}(\xi))$$ denotes the {\em additive canonical character} of $\F_q$ and 
	\begin{equation}\label{Fdef} F_{a_1,\ldots,a_s}(X)=\delta \sum_{i=1}^sa_if(X+\alpha_i).
	\end{equation}
	
	If $F_{a_1,\ldots,a_s}(X)$ is not of the form $g(X)^p-g(X)+c$ for any $(a_1,\ldots,a_s)\ne (0,\ldots,0)$,
	then the Weil bound, Lemma~\ref{weil}, can be applied and yields $(\ref{main})$.

	\section{Monomials $f_d(X)=X^d$}\label{monomials}
	Now we study the special case 
	$$f(X)=f_{dp^j}(X)=X^{dp^j} \quad \mbox{with}\quad \gcd(d,p)=1 \mbox{ and }j=0,1,\ldots$$
	
	Put ${\cal A}^k=\{\alpha^k: \alpha\in {\cal A}\}$.
	Since $(X+\alpha)^{dp^j}=(X^{p^j}+\alpha^{p^j})^d$ and $\xi\mapsto \xi^{p^j}$ permutes $\F_q$ we have
	$${\cal T}(\underline{c},{\cal A},f_{dp^j})={\cal T}(\underline{c},{\cal A}^{p^j},f_d)$$
	and we may assume $j=0$.
	Since 
	$$\xi^q=\xi\quad \mbox{for all }\xi\in \F_q$$ 
	we may restrict ourselves to the case $d<q$.

	To prove the first part of Theorem \ref{Xd} we have to show that $(\ref{main})$ is applicable. By $(\ref{Fdef})$ with 
	$$f(X)=f_d(X)=X^d$$ 
	we have
	$$F_{a_1,\ldots,a_s}(X)=\delta \sum_{i=1}^sa_i(X+\alpha_i)^d$$
	and thus
	\begin{equation}\label{deriv} F'_{a_1,\ldots,a_s}(X)=\delta d \sum_{\ell=0}^{d-1} {d-1\choose \ell} \left(\sum_{i=1}^sa_i\alpha_i^\ell\right)X^{d-\ell-1}.
	\end{equation}
	Assume that for some $(a_1,\ldots,a_s)\in \F_p^s\setminus\{(0,\ldots,0)\}$ we have 
	\[
	F_{a_1,\ldots,a_s}(X)=g(X)^p-g(X)+c
	\]
	for some polynomial $g(X)\in \F_q[X]$ and some constant $c\in \F_q$.
	We have
	\begin{equation}\label{dmodp} \mbox{either }F_{a_1,\ldots,a_s}(X)= c \quad \mbox{or}\quad
	\deg(F_{a_1,\ldots,a_s})\equiv 0\bmod p
	\end{equation}
	and
	$$F'_{a_1,\ldots,a_s}(X)=-g'(X).$$
	Then either
	\begin{equation}\label{F'1} F'_{a_1,\ldots,a_s}(X)= 0
	\end{equation}
	or
	\begin{equation}\label{F'2}
	\deg(F'_{a_1,\ldots,a_s})< \deg(g)= \frac{\deg(F_{a_1,\ldots,a_s})}{p}.
	\end{equation}
	Let
	\[
	d=d_0+d_1p+\dots+d_{r-1}p^{r-1},\quad 0\le d_0,\ldots,d_{r-1}<p,\quad d_0\not=0,
	\]
	be the $p$-adic expansion of $d$.
	Assume that there are $k\ge 0$ consecutive digits 
	\[
	d_m=d_{m+1}=\ldots=d_{m+k-1}=p-1,\quad 1\le m\le m+k\le r-1,
	\]
	of maximal size and
	\[
	s\le \max\left\{d_0,(d_{m+k}+1)p^k\right\}.
	\]
	Note that $\deg(F_{a_1,\ldots,a_s})\le d-d_0$ by $(\ref{dmodp})$ with the convention $\deg(0)=-1$. In both cases, $(\ref{F'1})$ and $(\ref{F'2})$, the coefficients of $F'_{a_1,\ldots,a_s}(X)$ at $X^{d-1-\ell}$ are zero for $\ell=0,\ldots,d-(d-d_0)/p-1$.
	Since $\delta d\not= 0$ we get from $(\ref{deriv})$
	\begin{equation}\label{withbinom}
	{d-1 \choose \ell}
	\left( \sum_{i=1}^s a_i\alpha_i^\ell\right)=0,\quad \ell=0,\ldots,d-(d-d_0)/p-1.
	\end{equation}
	By Lucas' congruence, Lemma~\ref{lm:lucas}, we have
	\begin{equation}\label{bin1}{d-1 \choose \ell}\equiv {d_0-1\choose \ell}\not\equiv 0\bmod p,\quad \ell=0,\ldots,d_0-1,
	\end{equation}
	as well as
	\begin{equation}\label{bin2}{d-1\choose p^m\ell}\not\equiv 0\bmod p,\quad \ell=0,\ldots,(d_{m+k}+1)p^k-1,  
	\end{equation}
	since
	$$
	d-1= e_0+(p-1)(p^m+\ldots+p^{m+k-1})+d_{k+m}p^{k+m}+e_1p^{k+m+1}$$
	for some 
	$$0\le e_0<p^m,~0\le e_1<p^{r-k-m-1},$$
	and
	$$p^m\ell = \ell_0p^m+\ldots+\ell_{k-1}p^{m+k-1}+\ell_kp^{m+k}$$
	for some
	$$0\le \ell_0,\ldots,\ell_{k-1}<p,~0\le \ell_k\le d_{k+m},$$
	and any $0\le \ell \le (d_{m+k}+1)p^k-1$.
	Combining $(\ref{withbinom})$ with $(\ref{bin1})$ and $(\ref{bin2})$, respectively, we get
	\begin{equation}\label{eqsyst1}
	\sum_{i=1}^s a_i\alpha_i^{\ell}=0, \quad \ell=0,\dots d_0-1,
	\end{equation}
	and 
	\begin{equation}\label{eqsyst2}
	\sum_{i=1}^s a_i\alpha_i^{p^m\ell}=0, \quad \ell=0,\ldots,(d_{m+k}+1)p^k-1,
	\end{equation}
	respectively. 
	Hence, if $s\le d_0$ or  $s\le (d_{m+k}+1)p^k$,
	the $s \times s$ coefficient matrix of the equations for $\ell=0,\ldots,s-1$ of $(\ref{eqsyst1})$ or $(\ref{eqsyst2})$, respectively, is an invertible Vandermonde matrix and we get 
	$$a_i=0,\quad i=1,\ldots,s,$$ contradicting $(a_1,\ldots,a_s)\in \F_p^s\setminus\{(0,\ldots,0)\}$. For the second case we used that $\xi\mapsto \xi^{p^m}$ permutes $\F_q$ and the $\alpha_i^{p^m}$, $i=1,\ldots,s$, are pairwise distinct.

	Proof of the second part of Theorem \ref{Xd}: Now assume $d<p^n$ for some $n$ with $1\le n\le r$, that is, $d_n=\ldots=d_{r-1}=0$, and assume $(\ref{Dcond})$ and $(\ref{scond2})$.
	Let $D$ be the number of binomial coefficients ${d\choose \ell}$, $\ell=1,\ldots,d$, which are nonzero modulo $p$.
	By Lemma \ref{th:fine} we have
	$$D=(d_0+1)\cdots(d_{n-1}+1)-1.$$
	For any $\alpha\in \F_q$ the polynomial 
	$$(X+\alpha)^d-\alpha^d= \sum_{\ell=0}^{d-1} {d\choose \ell}\alpha^\ell X^{d-\ell}$$
	is in the vector space generated by the monomials $X^{d-\ell}$ with nonzero ${d\choose \ell}\bmod p$, $\ell=0,\ldots,d-1$, of dimension $D$.
	For $D<s\le q$ and any subset $\{\alpha_1,\ldots,\alpha_s\}$ of $\F_q$ there is a nontrivial linear combination
	$$\sum_{i=1}^s\rho_i\left((X+\alpha_i)^d-\alpha_i^d\right)= 0$$ of the zero polynomial with 
	$(\rho_1,\ldots,\rho_s)\in \F_q^s\setminus\{(0,\ldots,0)\}$.
	If $D<s\le p$ and we take $\alpha_i\in \F_p$, $i=1,\ldots,s$, then we may assume 
	$\rho_i=a_i\in \F_p$ and 
	$$\sum_{i=1}^s a_i {\rm Tr}\left(\delta\left((\xi+\alpha_i)^d-\alpha_i^d\right)\right)=0\quad \mbox{for all }\xi\in \F_q,$$ 
	Taking $(a_1,\ldots,a_s)\in \F_p^s\setminus\{(0,\ldots,0)\}$ 
	from the previous step, the vector space of solutions $(c_1,\ldots,c_s)\in \F_p^s$ of the equation 
	$$a_1c_1+\ldots+a_sc_s=0$$
	is of dimension $s-1$.
	 More precisely, the mapping 
		$$\varphi:\F_p^ s \rightarrow \F_p, \quad \varphi(c_1,\ldots,c_s)=a_1c_1+\ldots+a_sc_s$$ 
		is surjective since $(a_1,\ldots,a_s)$ is not the zero vector. By the rank-nullity theorem its kernel is of dimension $s-1$.
	That is, not all $\underline{c}=(c_1,\ldots,c_s)\in \F_p^s$ are attained as $\underline{c}=\left({\rm Tr}\left(\delta\left((\xi+\alpha_i)^d-\alpha_i^d\right)\right)\right)_{i=1}^s$ for any $\xi\in \F_q$, 
	namely those $\underline{c}$ which are not in the kernel of $\varphi$.
	We can extend this argument to $s>p$ by extending $(a_1,\ldots,a_p)\in \F_p^p\setminus \{(0,\ldots,0)\}$ to $(a_1,\ldots,a_p,0,\ldots,0)\in \F_p^s\setminus\{(0,\ldots,0)\}$.

	Proof of the third part of Theorem \ref{Xd}: Now we drop the condition $(\ref{Dcond})$ but then $s$ has to satisfy the stronger condition $(\ref{Dcond2})$ instead of $(\ref{scond2})$.
	For each $\alpha \in \F_q$ we have 
	\[
	{\rm Tr}(\delta((X+\alpha)^d-\alpha^d))=\sum_{\ell=0}^{d-1} {d\choose \ell} {\rm Tr}(\delta\alpha^\ell X^{d-\ell}),
	\]
	since the trace is $\F_p$-linear, and thus it
	lies in the $\F_p$-linear space generated by the polynomials ${\rm Tr}(\beta_i X^{d-\ell})$ with nonzero ${d\choose \ell}$ modulo $p$, $i=1,\ldots,r$, $\ell=1,\ldots,d$, of dimension at most $Dr$, where $\{\beta_1,\ldots,\beta_r\}$ is a basis of~$\F_q$ over $\F_p$.
	Now let $s> Dr$, then for any subset ${\cal A}=\left\{\alpha_1,\dots \alpha_s \right\}$ of $\F_q$ consider the set of polynomials 
	\[
	\left\{ {\rm Tr}(\delta((X+\alpha_i)^d-\alpha_i^d)): i=1,\ldots,s \right\}.
	\]
	Since $s>Dr$ there is a nontrivial $\F_p$-linear combination
	\[
	\sum_{i=1}^{s}a_i{\rm Tr}(\delta((X+\alpha_i)^d-\alpha_i^d)) = 0
	\]
	of the zero polynomial.
	Now consider the linear subspace of solutions $(c_1,\ldots,c_s)\in \F_p^s$ of the equation $a_1c_1+\dots+a_sc_s=0$ which is of dimension~$s-1$. Let $\underline{c} \in \F_p^s$ be a point which does not lie in this linear subspace, then $\underline{c}$ is not attained as $\underline{c}=({\rm Tr}(\delta((\xi+\alpha_i)^d-\alpha_i^d)))_{i=1}^s$ for any $\xi\in \F_q$.

	\section{Rational functions $f_{-d}(X)=X^{-d}$}\label{sec:rational}

	Let $f_{q-d-1}(X)=X^{q-1-d}$ be a monomial of degree $q-d-1$, where $1\le d< q$. With the convention $0^{-1}=0$ we can identify $f_{q-d-1}(X)$ with the rational function $f_{-d}(X)=X^{-d}$. 
	Let $\gcd(d,q)=p^j$.
	Since 
	$$(X+\alpha)^{-p^j}=(X^{p^j}+\alpha^{p^j})^{-1}$$ 
	and $\xi\mapsto \xi^{p^j}$ permutes $\F_q$ we have
	$${\cal T}(\underline{c},{\cal A},f_{-dp^j})={\cal T}(\underline{c},{\cal A}^{p^j},f_{-d})$$
	and may restrict ourselves to the case $\gcd(d,q)=1$, that is,
	$$d=d_0+t_1p, \quad \mbox{where }1 \le d_0<p.$$
	
	We first show that there is no nonzero $s$-tuple 
	$$(a_1,\ldots,a_s)\in \F_p^s\setminus\{(0,\ldots,0)\}$$ 
	such that
	\[
	F_{a_1,\dots,a_s}(X)=\delta\sum_{i=1}^sa_i(X+\alpha_i)^{-d}=H(X)^p-H(X)
	\]
	for any rational function $H(X)\in \overline{\F_p}(X)$. We have
	\[
	F_{a_1,\dots,a_s}(X)=\frac{f(X)}{g(X)},
	\]
	where
	\[
	f(X)=\sum_{j=1}^sa_j\prod_{i\not=j}(X+\alpha_i)^d\]
	and
	\[g(X)=\prod_{i=1}^s(X+\alpha_i)^d.
	\]
	Suppose to the contrary that there exists a rational function 
	$$H(X)=\frac{u(X)}{v(X)}\in \overline{\F_p}(X)\quad \mbox{with }\gcd(u,v)=1\mbox{ and }v(X)\mbox{ is monic}$$
	satisfying
	$$F_{a_1,\ldots,a_s}(X)=H(X)^p-H(X).$$
	Therefore, we have
	\begin{equation}\label{ratform}
	\frac{f(X)}{g(X)}=\frac{u(X)^p}{v(X)^p}-\frac{u(X)}{v(X)}.
	\end{equation}
	Clearing denominators we obtain 
	$$f(X)v(X)^p=(u(X)^p-u(X)v(X)^{p-1})g(X)$$ 
	and thus $v(X)^p$ divides $g(X)$, hence 
	$$v(X)=\prod_{i=1}^s (X+\alpha_i)^{e_i}\quad \mbox{for some }0\le e_i\le t_1,~ 
	i=1,\ldots,s.$$ 
	Now by taking derivatives of both sides of \eqref{ratform} and clearing denominators we get
	\begin{equation}\label{deriveq}
	(f'(X)g(X)-f(X)g'(X))v(X)^2=(u(X)v'(X)-u'(X)v(X))g(X)^2.
	\end{equation}
	Without loss of generality we may assume $a_1\ne 0$, thus 
	$$f(-\alpha_1)=a_1\prod_{i=2}^s(\alpha_i-\alpha_1)^d\ne 0$$ 
	and 
	$$X+\alpha_1\mbox{ does not divide }f(X).$$ 
	Moreover, $(X+\alpha_1)^{d-1}$ and $(X+\alpha_1)^d$ are the largest powers dividing~$g'(X)$ and $g(X)$, respectively, that is, $$(X+\alpha_1)^{d-1+2e_1}$$ is the largest power of $(X+\alpha_1)$ dividing the left hand side of $(\ref{deriveq})$. Observing that $g(X)^2$ and thus the right hand side of $(\ref{deriveq})$ is divisible by 
	$$(X+\alpha_1)^{2d}$$ 
	we get $$d-1+2e_1\ge 2d$$ 
	and thus
	$$e_1\ge \frac{d+1}{2}> \frac{t_1p}{2}\ge t_1,$$ 
	which is a contradiction. 
	
	We showed that the conditions of Lemma~\ref{th:momo} are satisfied and Theorem~\ref{inv} follows
	from $(\ref{redchar})$ and Lemma~\ref{th:momo} since
	$$\left|\sum_{\xi\in \F_q}\psi(F_{a_1,\ldots,a_s}(\xi))\right|\le \left|\sum_{\xi\in \F_q\setminus -{\cal A}}\psi(F_{a_1,\ldots,a_s}(\xi))\right|+s,$$
	where $-{\cal A}=\{-\alpha_1,\ldots,-\alpha_s\}$.

	\section{Arbitrary polynomials}\label{secarb}
	In this section we prove Theorem~\ref{arb}.
	
	Let 
	$$f(X)=\sum_{j=0}^d \gamma_jX^j\in \F_q[X],\quad \gamma_d\ne 0,$$
	be a polynomial of degree 
	$$d=d_0+t_1p,\quad 1\le d_0<p,~0\le t_1<q/p.$$ 
	
	Proof of the first part: We have to show that $(\ref{main})$ is applicable, that is, the polynomial
	$F_{a_1,\ldots,a_s}(X)$ defined by $(\ref{Fdef})$ is not of the form $g(X)^p-g(X)+c$ for any~$(a_1,\ldots,a_s)\ne (0,\ldots,0)$.
	
	Suppose the contrary that there exists an $s$-tuple 
	$$(a_1,\ldots,a_s)\in \F_p^s\setminus\{(0,\ldots,0)\}$$
	such that the polynomial
	\[
	F_{a_1,\ldots,a_s}(X)=\delta\sum_{\ell=0}^d\left(\sum_{j=\ell}^d\sum_{i=1}^sa_i\gamma_j {j \choose \ell}\alpha_i^{j-\ell} \right)X^{\ell}
	\]
	can be written as 
	$$g(X)^p-g(X)+c\quad \mbox{for some }g(X)\in \F_q[X]\mbox{ and }c \in \F_q.$$ 
	We have either $$F_{a_1,\dots,a_s}(X)= 0$$ or
	\[
	\deg(F_{a_1,\ldots,a_s})\equiv 0\bmod p.
	\]
	Hence, 
	$$\deg(F_{a_1,\dots,a_s}) \le d-d_0,$$
	where we used the convention $\deg(0)=-1$.
	We conclude that the coefficients $\delta R_\ell$ of $F_{a_1,\dots,a_s}(X)$ at $X^{\ell}$ vanish for $\ell=d-d_0+1,\dots, d$.
	Since~$\delta\not= 0$ we have
	\begin{equation}\label{recu1}
	R_{\ell}=\sum_{j=\ell}^d\sum_{i=1}^sa_i\gamma_j {j \choose \ell}\alpha_i^{j-\ell}=0,\quad \ell=(d-d_0)+1,\ldots,d.
	\end{equation}
	Note that by Lucas' congruence, Lemma~\ref{lm:lucas},
	\begin{equation}\label{binr}{d\choose r}\equiv {d_0\choose r}\not\equiv 0\bmod p,\quad r=0,\ldots,d_0.
	\end{equation}
	Define $T_\ell$, $\ell=0,\ldots d_0-1$, recursively by 
	$$T_0=R_d$$ 
	and
	\begin{equation}\label{recu2}
	T_{\ell}=R_{d-\ell}-\gamma_d^{-1}\sum_{r=0}^{\ell-1}\gamma_{d-\ell+r} {r+d-\ell \choose d-\ell}{d \choose r}^{-1}T_r,
	\end{equation} 
	for $\ell=1, \ldots, d_0-1$.
	
	%
		Next we show that
		\begin{equation}\label{Tell0}
		T_{\ell}=\gamma_d{d \choose \ell}\sum_{i=1}^sa_i\alpha_i^\ell=0, \quad \ell= 0, \ldots, d_0-1.
		\end{equation}
		For $\ell=0$ the formula follows from \eqref{recu1}
		and for $\ell=1,\ldots,d_0-1$ from \eqref{recu2} we get by induction
		$$T_{\ell}=R_{d-\ell}-\sum_{r=0}^{\ell-1}\gamma_{d-\ell+r}{r+d-\ell\choose d-\ell} \sum_{i=1}^sa_i\alpha_i^r$$
		and from \eqref{recu1}
		$$T_\ell=\gamma_d{d\choose \ell}\sum_{i=1}^r a_i \alpha_i^\ell.$$
		Moreover, we get 
		$$T_\ell=0,\quad \ell=0,\ldots,d_0-1,$$
		from \eqref{recu1}, \eqref{recu2} again by induction.

	By \eqref{Tell0} and $(\ref{binr})$ we get since $\gamma_d\ne 0$,
	\[
	\sum_{i=1}^sa_i\alpha_i^{\ell}=0,\quad \ell=0,\ldots,d_0-1.
	\]
	Thus for $s\le d_0$, the $(s\times s)$-coefficient matrix 
	$$\left(\alpha_i^\ell\right)_{i=1,\ldots,s,\ell=0,1,\ldots,s-1}$$
	of the system of the first~$s$ equations is a regular Vandermonde matrix and we get $(a_1,\ldots,a_s)=(0,\ldots,0)$, which is a contradiction.

	For the second part of Theorem~\ref{arb} we  assume $f(X)\in \F_p[X]$ and notice that for any $\alpha\in \F_q$ the element 
	$f(X+\alpha)-f(\alpha)$ is in the vector space generated by the monomials~$X^i$, $i=1,\ldots,d$, of
	dimension $d$. 
	
	If $d<s\le p$, we can choose any ${\cal A}\subset \F_p$. Then 
	$$f(X+\alpha_i)-f(\alpha_i), \quad i=1,\ldots,s,$$ 
	are linearly dependent over $\F_p$ as well as 
	$${\rm Tr}(\delta(f(X+\alpha_i)-f(\alpha_i))), \quad i=1,\ldots,s,$$
	that is,
	$$\sum_{i=1}^sa_i{\rm Tr}(\delta(f(X+\alpha_i)-f(\alpha_i)))=0$$
	for some $(a_1,\ldots,a_s)\in \F_p^s\setminus\{(0,\ldots,0)\}$
	and the result follows since not all $(c_1,\ldots,c_s)\in \F_p^s$ satisfy $a_1c_1+\ldots+a_sc_s=0$.
	
	If $d<p$ and $s>p$, we can choose $(a_1,\dots, a_p) \in \F_p^p \setminus\{(0,\ldots,0)\}$ as in the case $s=p$ and extend it to $(a_1,\dots, a_p,a_{p+1},\ldots,a_s)\in \F^s_p \setminus\{(0,\ldots,0)\}$  with $a_{p+1}=\dots=a_s=0$.
	
	Proof of the third part of Theorem \ref{arb}: Recall that $\{\beta_1,\ldots,\beta_r\}$ is a basis of $\F_q$ over $\F_p$.
	Each $f(X+\alpha)-f(\alpha)$ lies in the $\F_p$-vector space generated by 
	$$\beta_jX^i, \quad j=1,\ldots,r,~i=1,\ldots,d,$$ 
	of dimension $dr$. The dimension of the vector space generated by 
	$${\rm Tr}(\beta_jX^i)\quad j=1,\ldots,r,~i=1,\ldots,d,$$ is at most $dr$. If $q\ge s>dr$, there is a nontrivial linear combination
	$$\sum_{i=1}^sa_i {\rm Tr}(\delta(f(X+\alpha_i)-f(\alpha_i)))=0$$ for any ${\cal A}\subset \F_q$ of size $s$ and the result follows.

	\section{Final Remarks}
	\label{final}
	
	\subsection{A counter-example for $\gcd(d,p)>1$}
	\label{counterp}
	
	Now we provide an example that the restriction $\gcd(d,q)=1$ in Theorem~\ref{arb} is needed.
	
	Choose any $f(X)$ of the form  
	$$f(X)=\delta^{-1}(g(X)^p-g(X)+c)\quad \mbox{for some }g(X)\in \F_q[X]\mbox{ and }c\in \F_q,$$
	$a_1=1$ and $a_2=\ldots=a_s=0$. 
	Then we obtain
	\[
	F_{a_1,\ldots,a_s}(X)=\delta f(X+\alpha_1)=g(X+\alpha_1)^p-g(X+\alpha_1)+c
	\]
	and thus
	\begin{eqnarray*}
		T(f(\xi+\alpha_1))&=&{\rm Tr}(\delta f(\xi+\alpha_1))={\rm Tr}(F_{a_1,\ldots,a_s}(\xi))\\
		&=&{\rm Tr}(g(\xi+\alpha_1)^p-g(\xi+\alpha_1)+c)={\rm Tr}(c)
	\end{eqnarray*}
	for all $\xi\in \F_q$,
	that is,
	any vector $(c_1,\ldots,c_r)\in \F_p^r$ with $c_1\ne {\rm Tr}(c)$ is not attained as $(T(f(\xi+\alpha_i)))_{i=1}^s$.
	
	We conclude that for polynomials of degree $d$ with $\gcd(d,p)>1$, the bound of Theorem~\ref{arb} may not hold. 
	However, by Theorem~\ref{Xd}, for monomials the restriction $\gcd(d,p)=1$ is not needed.
	
	\subsection{Missing digits and subsets}
	
	For subsets ${\cal D}$ of $\F_p$,
	the closely related problem of estimating the number of $\xi\in \F_q$
	with $$f(\xi)\in\{d_1\beta_1+\ldots+d_r\beta_r :d_1,\ldots,d_r\in {\cal D}\}$$
	was studied in \cite{damasa,diessh,ga}, that is, $f(\xi)$ 'misses' the digits in $\F_p\setminus {\cal D}$.
	It is straightforward to extend these results combining our approach 
	with certain bounds on character sums
	to estimate the number of
	$\xi\in \F_q$ with 
	$$f(\xi+\alpha_i)\in \{d_1\beta_1+\ldots+d_r\beta_r: d_1,\ldots,d_r\in {\cal D}\},\quad i=1,\ldots,s.$$
	For example, for ${\cal D}=\{0,\ldots,t-1\}$ we can use the bound on exponential sums of \cite{ko}.
	
	Instead of restricting the set of digits we may restrict the set of $\xi$. That is, for a subset ${\cal S}$
	of $\F_q$ we are interested in the number of solutions~$\xi\in {\cal S}$
	of 
	$$(T(f(\xi+\alpha_1)),\ldots,T(f(\xi+\alpha_s)))=\underline{c}$$
	for any fixed $\underline{c}\in \F_p^s$.
	Typical choices of ${\cal S}$ are 'boxes' \cite{dale,ko} and 'consecutive' elements \cite{wi}.

	\subsection{Optimality and prescribed digits}
	
	Swaenepoel~\cite{sw1} improved the bound \eqref{DarSar} of~\cite{DS13} in the case when the polynomial $f(X)$ has degree~$2$ or is a monomial. In particular, for $s=1$ and $d=2$ the improved bound of \cite{sw1} is optimal.
	She also generalized \eqref{DarSar} to several polynomials with $\F_p$-linearly 
	independent leading coefficients \cite[Theorem~5.1]{sw1}.
	
	Moreover, in \cite{sw2} Swaenepoel studied the number of solutions $\xi\in \F_q$ for which some of the digits of $f(\xi)$ are prescribed, that is, for given ${\cal I}\subset\{1,\ldots,r\}$ and given $c_i\in \F_p$, $i\in {\cal I}$, the number of $\xi\in \F_q$
	with
	$${\rm Tr}(\delta_if(\xi))=c_i, \quad i\in {\cal I}.$$

	\subsection{Related work on pseudorandom number generators}
	Some of the ideas of the proofs in this paper are based on earlier work on nonlinear, in particular, inversive pseudorandom number generators, see~\cite{MW03,niwi00,niwi01}.
	
	More precisely, in \cite{niwi01} the $q$-periodic sequence $(\eta_n)$ over $\F_q$ defined by 
	$$\eta_{n_1+n_2p+\dots+n_rp^{r-1}}=f(n_1\beta_1+\ldots+n_r\beta_r),\quad 0\le n_1,\ldots,n_r<p,$$
	passes the {\em $s$-dimensional lattice test} if $s$
	polynomials of the from 
	$$f(X+\alpha_j)-f(X),\quad j=1,\ldots,s,$$
	are $\F_q$-linearly independent. However, in the proofs of this paper we need that they are linearly independent over $\F_p$.
	
	To prove Theorem~\ref{inv} for $d<p$, the method of \cite{niwi00} can be easily adjusted using \cite[Lemma~2]{niwi00}. However, for $d\ge p$ we had to use a different approach since \cite[Lemma~2]{niwi00} is not applicable in this case.
	
	Finally, in the proof of \cite[Theorem~4]{MW03} we showed that polynomials of the form $F_{a_1,\ldots,a_s}(X)$ can only be identical $0$ if $a_1=\ldots=a_s=0$. However, in the proof of Theorem~\ref{arb} we had to show that
	$F_{a_1,\ldots,a_s}(X)$ is not of the form $g(X)^p-g(X)+c$
	and we had to modify the idea of \cite{MW03}.

	\subsection{Rudin-Shapiro function}
	
	The {\em Rudin-Shapiro sequence} $(r_n)$ is defined by 
	$$r_n=\sum_{i=0}^\infty n_in_{i+1},\quad n=0,1,\ldots$$
	if 
	$$n=\sum_{i=0}^\infty n_i2^i,\quad n_0,n_1,\ldots\in \{0,1\}.$$
	M\"ullner showed that the Rudin-Shapiro sequence along squares $(r_{n^2})$ is normal \cite{mu}.
	
	The \emph{Rudin-Shapiro function} $R(\xi)$ for the finite field $\F_q$ with respect to the ordered basis $(\beta_1,\ldots,\beta_r)$ is defined as 
	\[
	R(\xi)=\sum_{i=1}^{r-1}x_ix_{i+1},\quad \xi=x_1\beta_1+x_2\beta_2+\dots+x_r\beta_r,\quad x_1,\ldots,x_r\in \F_p.
	\]
	For $f(X)\in \F_q[X]$ and $c\in \F_p$ let 
	\[
	{\cal R}(c,f)=\left\{\xi \in \F_q: R(f(\xi))=c \right\}.
	\]
	It seems to be not possible to use character sums to estimate the size of~${\cal R}(c,f)$.
	However, in \cite{damewi} the Hooley-Katz Theorem, see \cite[Theorem~7.1.14]{mupa} or \cite{hoka} was used to show that if $d=\deg(f)\ge 1$, 
	\[
	\left|{\cal R}(c,f)-p^{r-1}\right|\le C_{r,d}p^{\frac{3r+1}{4}},
	\]
	where $C_{r,d}$ is a constant depending only on $r$ and $d$.
	In particular, we have for fixed $d$ and $r\ge 6$,
	$$\lim_{p\rightarrow\infty}\frac{|{\cal R}(c,f)|}{p^{r-1}}=1,$$
	that is, $R(f)$ is $p$-normal for $s=1$ and $r\ge 6$.
	
	However, we are not aware of a result on the $r$-normality of $R(f)$.

\section*{Acknowledgment}
The authors are partially supported by the Austrian Science Fund FWF Project P 30405.
They wish to thank the anonymous referee for very useful suggestions.

\end{document}